\documentclass[12pt]{article}

\newcommand{\floor}[1]{\left\lfloor #1 \right\rfloor}

\renewcommand{\r}[2][p]{r_{#1}\left(#2\right)}

\newcommand{\fq}[1]{\left\lfloor Q_{#1} \right\rfloor}
\newcommand{\fqm}{\left\lfloor Q_m \right\rfloor}

\newcommand{\fqM}{\left\lfloor Q_M \right\rfloor}

\newcommand{\fqMm}{\left\lfloor Q_{M-1} \right\rfloor}

\newcommand{\fr}[1]{\left\lfloor R_{#1} \right\rfloor}
\newcommand{\frm}{\left\lfloor R_m \right\rfloor}

\newcommand{\frM}{\left\lfloor R_M \right\rfloor}

\newcommand{\dqr}[1][m]{\Big(\left\lfloor Q_{#1} \right\rfloor-\left\lfloor R_{#1} \right\rfloor\Big)}
\newcommand{\dqrsp}[1][m]{\big(\lfloor Q_{#1} \rfloor-\lfloor R_{#1} \rfloor\big)}
\newcommand{\drqsp}[1][m]{\big(\lfloor R_{#1+1} \rfloor-\lfloor R_{#1} \rfloor\big)}
\newcommand{\drq}[1][m]{\Big(\left\lfloor R_{#1+1} \right\rfloor-\left\lfloor Q_{#1} \right\rfloor\Big)}
\newcommand{\drr}[1][m]{\Big(\left\lfloor R_{#1+1} \right\rfloor-\left\lfloor R_{#1} \right\rfloor\Big)}

\newlength{\largo}
\newlength{\corto}
\newlength{\resta}

\newcommand{\inviss}[2]{ &=\ds{#1} \\
                         &
                  \settowidth{\largo}{$\ds{=#1}$} \settowidth{\corto}{$\ds{#2}$} \setlength{\resta}{\largo}\addtolength{\resta}{-\corto}
                  \,\,\hspace*{\resta}\ds{#2} }

\newcommand{\quadf}{\Q(\sqrt{-p})}
\usepackage{float}
\restylefloat{table}

\newcommand {\ds}{\displaystyle}

\newtheorem{thm}{Theorem}[section]
\newtheorem{lem}[thm]{Lemma}
\newtheorem{cor}[thm]{Corollary}

\newtheorem{rem}[thm]{Remark}
\newtheorem{nota}[thm]{Notation}
\newtheorem{eg}[thm]{Example}

\newtheorem{dfn}[thm]{Definition}

\newcommand{\ok}{\hfill $\Box$\\[2mm]}
\newcommand{\okok}{\hfill \Box}
\newlength{\tamano}
\newlength{\taman}
\newcounter{figura}
\newcommand{\crea}[2]{\parbox[l][#2cm][l]{#1cm}{\ }}

\newcommand{\at}[2][rrrrrrrrrrrrrr]{    
                     \begin{array}{#1}
                     #2\\
                     \end{array}}
\newcommand{\defi}[2]{#1 = \left\{\at[lllllllllllll]{#2}\right.} 

\newcommand{\sx}[1][.75]{\crea{0}{#1}}

\def\newsqrt{\mathpalette\DHLhksqrt} 
\def\DHLhksqrt#1#2{\setbox0=\hbox{$#1\sqrt{#2\,}$}\dimen0=\ht0
\advance\dimen0-0.45\ht0
\setbox2=\hbox{\hspace{-.014cm}\vrule height\ht0 depth -\dimen0}%
{\box0\lower0.55pt\box2}\,}

\newcommand{\R}{{\mathbb R}}

\newcommand{\Q}{{\mathbb Q}}
\newcommand{\N}{{\mathbb N}}

\newcommand{\Z}{{\mathbb Z}}

\usepackage{amsfonts,amsmath,enumerate}

 \usepackage{graphicx}
 \usepackage{color}
 \usepackage{enumerate}
 \usepackage{epsfig}
 \usepackage{fancyhdr}
 \usepackage{ifthen}
 \usepackage{mathtools}
 \usepackage[tiling]{pst-fill}
 \usepackage{xifthen}
 \usepackage{xparse}
 \usepackage{tikz}
 \usepackage{float} 

\begin{document}
\title{A computable formula for the class number of the imaginary quadratic field $\quadf, \ p=4n-1$}
\author{Jorge Garcia}
\maketitle


\begin{abstract}
   A formula for the class number $h$ of the imaginary quadratic field $\quadf$ is obtained by counting on a specific way the quadratic residues of a prime number of the form $p=4n-1.$
   Formulas for the sum of the quadratic residues are also found.
\end{abstract}



\section{Introduction}\label{Section:Intro}
Given a prime number $p$ and $0<k<p$, we look at the residues of $k^2$ modulus $p.$  When we add those residues we obtain the
sum of quadratic residues relative to the prime number $p$. The question is, how does that sum behave? Is there a formula for such sum? What is that sum related to?
Here are the answers to these simple questions. There is indeed a very simple formula for this sum when the number is of the form $p=4n+1,$ in fact this sum is $\binom{p}{2}$. There is also a formula when the prime is of the form $4n-1,$ in fact this sum is $\binom{p}{2}-p\cdot h$ where $h$ is the class number of the imaginary quadratic field $\quadf$. Formulas for the class number when the prime is of the form $p=4n-1>3$  are

     \begin{eqnarray*}
       h &=& \frac{\sqrt{p}}{2\pi}\sum_{r=1}^{\infty}\frac{\chi(r)}{r} \\
         &=& \frac{-1}{4p}\sum_{r=1}^{4p}r\cdot\Big(\frac{p}{r}\Big),
     \end{eqnarray*}
where $\Big(\frac{p}{r}\Big)$ is the Kronecker symbol and $\chi$ is the Dirichlet character.
The first formula is the Dirichlet Class Number Formula and can be found in~\cite{Dirichlet:bedesat37} and the second formula can be found in \cite{Narkiewicz:eatan04} (Corollary 1 p. 428).

In Section \ref{Section:BasicLemmas} we will study the basic lemmas as well as the definition whereas in Section \ref{Section:MainResult} we will establish the main result for the sum of quadratic residues. In Section \ref{Section:Computable} we establish several formulas for the class number when $p$ is of the form $4n-1$ and finally in Section \ref{Section:Examples} we provide some examples and conjecture some results with a mystery term.


\section{Sum of Quadratic Residues}\label{Section:BasicLemmas}
If we look at small primes of the form $p=4n-1$ and we add the residues of $k^2$ and $(2n-k)^2$ mod $p$ we observe some patterns as $k$ ranges from 1 to $2n$. Unfortunately these patterns do not hold for larger primes, however, we need to make some adjustments to those patterns. The following lemmas will capture those patterns and these lemmas will be the basic stones of our main theorems.

\begin{dfn}\label{residues}
  Let $q$ be a positive integer and $x\in\Z$. By $r_q(x)$ we denote the residue of $x$ when we divide by $q$. Hence $r_q(x)\in\{0,1,2,...,q-1\}$ satisfies
  \[x=h\cdot q + r_q(x),\]
  for some $h\in \Z$. Clearly, $h=\floor{x/q}.$
\end{dfn}

\begin{lem}\label{Lemma1Patternsp=4n-1}
Let $n\in\N,\  p,h,r,k\in\Z$ such that $p$ is prime, $p=4n-1, \ k^2=hp+r, \ 0\le r\le 4n-2$ and $0\le k\le n.$
\begin{enumerate}[(i)\ ]
  \item If $k\ge r-3n+2$ then
       \[ \r{k^2}+\r{(2n-k)^2}=2r-k+n. \]
  \item If $k<r-3n+2$ then
       \[ \r{k^2}+\r{(2n-k)^2}=2r-k-3n+1. \]
\end{enumerate}
\end{lem}

\noindent {\bf Proof.}
  Notice that
    \[
       \begin{split}
          (2n-k)^2 & \equiv 4n^2-4nk+k^2 -(4n-1)(n-k) \pmod{p} \\
                   & = k^2+n-k = hp+r+n-k \\
                   & \equiv r+n-k \pmod{p}.
       \end{split}
    \]
    Clearly $0\le r+n-k.$ \  If $k\ge r-3n+2,$ then \ $r+n-k\le 4n-2,$ hence
       \[ \r{k^2}+\r{(2n-k)^2}=2r-k+n. \]
    This proves the first part. For the second part notice also that
    \[
       \begin{split}
          (2n-k)^2 & \equiv  r+n-k -4n+1 \pmod{p}\\
                   & = r-k-3n+1 \pmod{p}.
       \end{split}
    \]
    Clearly $r-k-3n+1\le 4n-2.$ \  Since $k< r-3n+2,\  r-k-3n+1=r+n-k-4n+1>-1,$ hence
       \[ \r{k^2}+\r{(2n-k)^2}=2r-k-3n+1. \]  \ok

\begin{lem}\label{Theorem2LemmaSplits4n-1} 
Let $n\in\N,\  p, k, m\in\Z$ such that $p$ is prime, $p=4n-1$ and $0\le k\le n.$
\begin{enumerate}[(i)\ ]
  \item If \ $mp-k<k^2-k\le (m+1)p-n-1$ then
       \[ \r{k^2}+\r{(2n-k)^2}=2k^2-k+n-2mp. \]
  \item If \ $(m+1)p-n-1<k^2-k< (m+1)p-k$ then
       \[ \r{k^2}+\r{(2n-k)^2}=2k^2-k+n-(2m+1)p. \]
\end{enumerate}
\end{lem}

\noindent {\bf Proof.}
\begin{enumerate}[(i)\ ]
  \item Clearly \ $mp<k^2\le (m+1)p+k-n-1\le (m+1)p-1.$ Hence there is $r\in\{1,2,..,p-1\}$ such that $k^2=mp+r.$  Since
    \[
       \begin{split}
          r-3n+2 & =  k^2-mp-3n+2\\
                   & = k^2-(m+1)p+n+1 \le k.
       \end{split}
    \]
    By Lemma \ref{Lemma1Patternsp=4n-1},
       \[
          \begin{split}
             \r{k^2}+\r{(2n-k)^2} & = 2r+n-k=2(k^2-mp)+n-k\\
                  & =2k^2-k+n-2mp.
          \end{split}
        \]
  \item Clearly \ $mp<(m+1)p-n-1+k<k^2.$ \ Since $k^2<(m+1)p,\ mp<k^2<(m+1)p.$ \ Hence there is $r\in\{1,2,..,p-1\}$ such that $k^2=mp+r.$ \ Since $r-3n+2=k^2-mp-3n+2=k^2-(m+1)p+n+1>k,$ by Lemma \ref{Lemma1Patternsp=4n-1},
       \[
          \begin{split}
            \r{k^2}+\r{(2n-k)^2} & = 2r-k-3n+1 \\
                                 & = 2k^2-2mp-k-p+n\\
                                 & = 2k^2-k+n-(2m+1)p.\okok
          \end{split}
       \]
\end{enumerate}

\begin{cor}\label{Corollary3SquareRoots}
Let $n\in\N,\  p, k,\in\Z$ such that $p$ is prime, $p=4n-1$ and $0\le k\le n.$ Consider an integer $m\ge 0$.
\begin{enumerate}[(i)\ ]
  \item If \ $\newsqrt{mp}<k\le \frac12+\frac12\newsqrt{1+4\left[(m+1)p-n-1\right]}$ then
       \[ \r{k^2}+\r{(2n-k)^2}=2k^2-k+n-2mp. \]
  \item If \ $\frac12+\frac12\newsqrt{1+4\left[(m+1)p-n-1\right]}<k<\newsqrt{(m+1)p}$ then
       \[ \r{k^2}+\r{(2n-k)^2}=2k^2-k+n-(2m+1)p. \]
\end{enumerate}

\end{cor}
\noindent {\bf Proof.}  Notice that $\sqrt{m}<k$ \ implies \  $mp-k<k^2-k$, \ also \ $k\le \frac12+\frac12\newsqrt{1+4\left[(m+1)p-n-1\right]}$ implies $4k^2-4k\le 4\left[(m+1)p-n-1\right]$, hence
    \[mp-k<k^2-k\le(m+1)p-n-1\]
and by Lemma \ref{Theorem2LemmaSplits4n-1} we have the result. The second part is done similarly.
\ok

\bigskip

\noindent The following notation will be very useful for our future theorems and identities.

\begin{nota}\label{NotationQmRmM} For $n\in\N,\  m,\in\Z$ with $p$ is prime, $p=4n-1$ and $m\ge 0$ we denote
  \[
          \begin{split}
            Q_m & = \frac12+\frac12\newsqrt{1+4\left[(m+1)p-n-1\right]} = \frac12+\frac12\newsqrt{4mp+3p-4} \\
                               R_m  & = \newsqrt{mp}\crea{0}{1.25}\\
                               M  & = \floor{\frac{n^2}{p}}.
          \end{split}
       \]

\end{nota}

\noindent The following lemma will allow us to estimate the quadratic residues.

\begin{lem}\label{Lemma4Estimater}
Using the previous notation, consider $n\in\N, \ p=4n-1$ a prime number  and $r\in\{0,1,...,p-1\}$ such that \[n^2=Mp+r.\]
Assume that \ $M>0.$  Then \ $r\ge 1$. If $n=11$ then $r=4n-9$ and if $n\neq 11$ then $r\le 4n-10.$
\end{lem}
\noindent {\bf Proof.}
   Clearly $r\ge 1$ otherwise $p$ divides $n$ which is impossible. If $n=11$, then $p=43$ and hence \  $n^2=121=2\cdot 43+35$, \ therefore \ $r=35=4\cdot n-9.$ \ \ Assume now that \ $n\neq 11.$ \ Let \ $y=M+1$. \ Consider \ $u=4n-r,\ \ u\ge 2$ \ and \ $f(y)=4y^2-(y+u-1).$

   \begin{enumerate}[\underline{Case}]
     \item \underline{$u=3$.} \ Then $n^2=Mp+4n-3$. Hence \ $(n-3)(n-1)=n^2-4n-3=Mp$, \ therefore $p$ divides $(n-3)$ or $p$ divides $(n-1).$  This forces $n=3, p= 11, M=\floor{\frac{9}{11}}=0$ \ or \  $n=1, p=3, M=0,$ which is impossible.
     \item \underline{$u=4$.} \  Then $(n-2)^2=Mp$. This forces $n=2, p= 7, M=\floor{\frac{4}{7}}=0,$ \ which is also impossible.
   \end{enumerate}
   Assume now that $u\neq 3,4.$\  Now
    \begin{equation}\label{Equationny}
       n^2=Mp+4n-u=(M+1)p+1-u=y(4n-1)+1-u,
    \end{equation}
   hence \ $n$ \ satisfies \ $n^2-n\cdot 4y+(y+u-1)=0$,\  therefore
   \[
      n=2y\pm \sqrt{4y^2-(y+u-1)},
   \]
   hence there is a non-negative integer $x$ such that $f(y)=x^2.$  \ \ If $x=0$ then \ $n=2y$ \  and \ $0=4y^2-(y+u-1)$. From Equation~\ref{Equationny}, we have $4y^2=y(4n-1)+1-u$, \ hence $y(4n-1)=y$. This forces $y=0$ and  $n=0$ contradicting $n\ge 1.$.  Hence necessarily \ $x>0.$ Therefore \ $4y^2-x^2=y+u-1\ge 1$\ and then \ $y+x\le u-1$.  Hence $y\in\{0,1,2,...,u-2\}.$ \ We know, $y\neq 0$. If $y=1$ then $M=0$ contrary to our assumption. Therefore  $y\in \{2,...,u-2\}.$   We now continue analyzing the cases $u=2,5,6,7,...$

   \begin{enumerate}[\underline{Case}]
     \item \underline{$u=2$.} \ Then \ $y+x\le 1$, this is impossible as $x\ge 1$ and $y\ge 2.$
     \item \underline{$u=5$.} \ Then \ $y\in\{2, 3\}.$ \  Clearly $f(2)=10\neq x^2$ \ and \  $f(3)=29\neq x^2.$
     \item \underline{$u=6$.} \ Then \ $y\in\{2, 3, 4\}.$ \  If $y=2$, then $f(2)=9= x^2,$ \ this forces \  $x=3, n=1, p=3, M=0$ \ which is impossible.  Clearly $f(3)=28\neq x^2$ \ and \  $f(4)=55\neq x^2.$
     \item \underline{$u=7$.} \ Then \ $y\in\{2, 3, 4, 5\}.$ \  Clearly $f(2)=8, f(3)=27, f(4)=54$ and $ f(5)=89$\ none of which is a perfect square.
     \item \underline{$u=8$.} \ Then \ $y\in\{2, 3, 4, 5, 6\}.$ \  Clearly $f(2)=7, f(3)=26, f(4)=53, f(5)=88$ and $ f(6)=131$\ none of which is a perfect square.
     \item \underline{$u=9$.} \ Then \ $y\in\{2, 3, 4, 5, 6, 7\}.$ \  Clearly $f(2)=6, f(4)=52, f(5)=87, f(6)=120$ and $ f(7)=181$\ none of which is a perfect square. This leaves us with the case $y=3$, which forces $f(3)=25, x=5, n\in\{11,1\}.$ But we know $n\neq 11$, hence $n=1, M=0,$ which is impossible.
   \end{enumerate}
   Therefore $u\ge 10$ and hence $r\le 4n-10.$
\ok

The following lemma provides the computation of quadratic residues when $n$ is congruent to 0,1,2 and 3 modulus 4 and $p=4n-1.$ Such lemma will be useful when we provide some formulas regarding quadratic residues.
\begin{lem}\label{Lemma4.5Chart}
  Let $n, M, p$ and $r$ as in Lemma~\ref{Lemma4Estimater}. Then $M=\floor{\frac{n}{4}}$ \ and the quadratic residue $r$ is given in Table~\ref{TableResiduesn} for each possible $n$ of the form $4k, 4k+1, 4k+2$ and $4k+3.$
 \begin{table}[H]
 \centering
 \begin{tabular}{|c|c|c|c|c|}
   \hline
   $n$ & $M$ & $p$ & $r$ \\   \hline  \hline
    $4k$& $k$ & $16k-1$ & $M$ \\  \hline
    $4k+1$ & $k$ & $16k+3$ & $5M+1$ \\  \hline
    $4k+2$ & $k$ & $16k+7$ & $9M+4$ \\  \hline
    $4k+3$ & $k$ & $16k+11$ & $13M+9$ \\  \hline
 \end{tabular}
 \caption{The quadratic residues $r=\r{n^2}$ for each of the four different cases of $n$.}\label{TableResiduesn}
 \end{table}

\end{lem}
\noindent {\bf Proof.}
We will verify that $\ds{M=\floor{\frac{n^2}{4n-1}}=k}$ \ in one case,  the other cases are done similarly. Let  $n=4k, \ k\ge 1.$ \ Clearly
     \[
         k\cdot \big(4(4k)-1\big) \le 16k^2 < 16k^2+15k-1=(k+1)\cdot \big(4(4k)-1\big),
     \]
hence $\ds{k\le \frac{n^2}{4n-1}< k+1}$, therefore $\ds{M=\floor{\frac{n^2}{4n-1}}=k}$. Now, observe that
       \[
          \begin{split}
              (4k)^2& = k\cdot \big(4(4k)-1\big) +k,        \\
            (4k+1)^2& = k\cdot \big(4(4k+1)-1\big) +5k+1,   \\
            (4k+2)^2& = k\cdot \big(4(4k+2)-1\big) +9k+4,   \\
            (4k+3)^2& = k\cdot \big(4(4k+3)-1\big) +13k+9.
          \end{split}
       \]
We can see that each of the numbers $r$ in Table~\ref{TableResiduesn} satisfies $n^2=Mp+r$ \ in each of the four different cases.  \ok

\bigskip

The following lemma will allow us to place $n-1, n, n+1$ in the right interval as established in Corollary~\ref{Corollary3SquareRoots}. This lemma and the following one will be useful to the proof of our main theorem which will come right after.

\begin{lem}\label{Lemma5QmRm}
  Let $n\in \N$ and $Q_m, R_m$ and $M$ as in Notation~\ref{NotationQmRmM} and $r=\r{n^2}$. Then
  \begin{enumerate}[(i)\ ]
    \item For $m\le M-1, \ \ R_m<Q_m<R_{m+1}.$
    \item If $M\ge 1$ \  then \ $R_{M-1}<n-1<Q_M.$
    \item If $M\ge 1$ \  then \ $R_{M-1}<n-2<Q_{M-1}.$
    \item $R_{M}<n<Q_M.$
    \item $R_{M}<n+1<Q_{M+1}.$
    \item $r\neq 2n-3,$ \ and \ $r\neq 2n-2.$ \ If $n>1$ then $r\neq 2n-1.$
    \item If $M\ge 1$ \  then \ $n-1\not\in(Q_{M-1},R_M).$
    \item If $M\ge 1$ \  then \ $n-1<Q_{M-1} \implies  n+1<Q_M \implies n+1\not\in(Q_M,R_{M+1})$.\\ Also $n-1<Q_{M-1} \Longleftrightarrow r<2n-3.$
    \item $R_M<n-1 \implies R_{M+1} < n+1 \implies n+1\not\in(Q_M,R_{M+1}).$
    \item $Q_M<R_{M+1} \Longleftrightarrow 2n-2<r.$
    \item $Q_M<n+2, \  R_{M+1}<n+2.$
  \end{enumerate}
\end{lem}

\medskip

\noindent {\bf Proof.}

  \begin{enumerate}[(i)\ ]
    \item Notice that $R_m<Q_m \Longleftrightarrow 2<3n$ which is true for all $n$. Notice also that $Q_m<R_{m+1} \Longleftrightarrow (m+1)(4n-1)<n^2+2n+1.$ Since $m\le M-1$, we have that $(m+1)(4n-1)\le Mp=n^2-r<n^2+2n+1.$
    \item Notice that $R_{M-1}<n-1 \Longleftrightarrow (M-1)p<n^2-2n+1.$ We certainly have
          $Mp-p=n^2-r-p<n^2-2n+1$. \ Also \ $n-1<Q_M \iff n^2-3n+4 < (M+1)p-n-1 \iff n^2+4<n^2-r+p+2n-5$,  and certainly by Lemma~\ref{Lemma4Estimater}, $r<6n-10.$
    \item We observe that $R_{M-1}<n-2 \iff (M-1)p<n^2-4n+4$. We do have $Mp-p<n^2-r-4n+1<n^2-4n+4.$ Likewise, $n-2<Q_{M-1}\iff n^2-5n+6<Mp-n-1,$ we certainly have, by Lemma~\ref{Lemma4Estimater}, $n^2-5n+6<n^2-r-n-1.$
    \item Observe that $R_{M}<n$ is equivalent to $Mp<n^2$ which is definitely true. Also $n<Q_M$ is equivalent to $n^2-n<Mp+p-n-1=n^2-r+p-n-1$ which is true as $r<p-1.$
    \item Note that $R_{M}<n+1$ is equivalent to $Mp<n^2+2n+1$ which is true because $n^2-r<n^2+2n+1$. Also \ $n+1<Q_M$ \ is equivalent to \ $n^2+n<(M+2)p-n-1=n^2-r+2p-n-1,$ this last inequality is true by Lemma~\ref{Lemma4Estimater}.
    \item \label{ImpossiblCaseLemma} If $r=2n-2$ then
        \begin{equation}\label{EquationCase6}
              n^2+n(-4M-2)+M+2=0,
        \end{equation}
        hence there is a non-negative integer $x$ such that $(2M+1)^2-(M+2)=x^2.$  Thus $2M+1+x\le (2M+1)^2-x^2= M+2,$ therefore $M+x\le 1.$ This gives us three possibilities $M=0,x=0$ or $M=0,x=1$ or $M=1,x=0.$ From $n=2M+1\pm x$ and Equation~\ref{EquationCase6}, each of the previous possibilities leads us to a contradiction.

        \indent If $r=2n-3$ then $Mp=n^2-2n+3$ hence $p$ divides $4n^2-8n+12.$ Since $p$ divides $4n^2-n$ and $8n-2$, we conclude that \ $p$ divides $n+10$, which is not possible. Now let $n>1,$ and assume $r=2n-1.$ Then $Mp=(n-1)^2$ which forces $M=0, n=1,$ contradicting our assumption.
    \item If $n\in(Q_{M-1},R_M),$ then $Mp-n-1<n^2-3n+2$ and $n^2-2n+1<Mp=n^2-r$, this forces $2n-3<r<2n-1$, i.e. $r=2n-2$, which was proven impossible in case \ref{ImpossiblCaseLemma}. Hence $n\in(Q_{M-1},R_M)$ is impossible.
    \item Clearly $n-1<Q_{M-1} \iff r<2n-3.$ \  Now, $n+1<Q_M$ iff\  $r\le 2n-3$  iff \  $r<2n-3$ as by case \ref{ImpossiblCaseLemma},  $r\neq 2n-3.$
    \item Observe that $R_M<n-1 \implies r>2n-1$ and $R_{M+1} < n+1 \iff Mp+p<n^2+2n+1 \iff 2n-2<r.$ The results follows.
    \item Note that $Q_M<R_{M+1} \iff (M+1)p<n^2+2n+1 \iff 2n-2<r.$
    \item Observe that $R_{M+1}<n+2$ is equivalent to $Mp<n^2+4n+5$ which is the same as $0<r+4n+5,$ of course this last inequality is true. Finally, $Q_M<n+2$ is equivalent to $Mp<n^2+4$ which is the same as the true inequality $0<r+4.$
  \end{enumerate}\ok

\begin{rem}\label{Remark6FlooryMinusFloorx}
  Let $x,y$ two real numbers. Then
  \[\big|\left\{k\in\Z\ :\ x<k\le y \right\}\big| =\floor{y}-\floor{x}.\]
\end{rem}

\begin{rem}\label{Remark7FlooryMinusFloorx}
  Let $x\in\R,\ y\in\R-\Z.$ \ Then
  \[\big|\left\{k\in\Z\ :\ x<k< y \right\}\big| =\floor{y}-\floor{x}.\]
\end{rem}

\bigskip

\begin{lem}\label{Lemma8EitherOr}  
 Let $n\in \N$ and $Q_m, R_m$ and $M\ge 1$ as in Notation~\ref{NotationQmRmM} and $r=\r{n^2}$. Then
 \begin{enumerate}[(i)\ ]
   \item Either \ $n-1=\floor{Q_{M-1}}=\floor{R_M}$ \ or \  $n-2=\floor{Q_{M-1}}=\floor{R_M}.$
   \item Either \ $n+1=\floor{Q_{M}}=\floor{R_{M+1}}$ \ or \  $n=\floor{Q_{M}}=\floor{R_{M+1}}.$
 \end{enumerate}
\end{lem}
\noindent {\bf Proof.}  From Lemma~\ref{Lemma5QmRm}, $n-1\in(R_{M-1},Q_M)-(Q_{M-1},R_M)=(R_{M-1},Q_{M-1}]\cup [R_M,Q_M).$

 \begin{itemize}
   \item  Assume now $n-1\in (R_{M-1},Q_{M-1}]$. \ If $n-1=Q_{M-1}$ \ then \ $r=2n-3$  \, which is not possible, hence $n-1<Q_M.$  By Lemma~\ref{Lemma5QmRm}, $n+1<Q_M$ \, and \, $R_M<n.$ Therefore \ $n-1=\floor{Q_{M-1}}=\floor{R_M}.$ The situation looks like the one on Figure~\ref{FigureFirstSituation}.

            \begin{figure}[H]
             \begin{center}
               \begin{tikzpicture}
                  \draw[thick] (0,0) -- (9,0);
                  \draw[thick] (4.5,.35) arc (160:200:1cm) node [below=5pt] {$R_M$};
                  \draw[thick] (2.5,.35) arc (20:-20:1cm) node [below=5pt] {$Q_{M-1}$};
                  \draw[thick] (7.5,.35) arc (20:-20:1cm) node [below=5pt] {$Q_M$};
                  \fill (1,0) circle (2.5pt) node [above=5pt] {$n-1$};
                  \fill (5,0) circle (2.5pt) node [above=5pt] {$n$};
                  \fill (6,0) circle (2.5pt) node [above=5pt] {$n+1$};
               \end{tikzpicture}
               \caption{The resulting case when $n-1\in (R_{M-1},Q_{M-1}]$.}{\label{FigureFirstSituation}}
             \end{center}
            \end{figure}
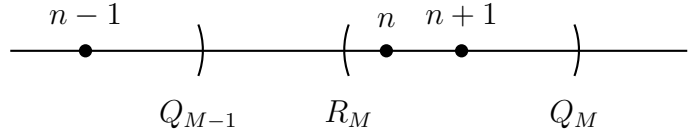
       If $2n-2<r$, by Lemma~\ref{Lemma5QmRm} (x),(xi), $Q_M<R_{M+1}<n+2.$ \ If $2n-2>r$, by Lemma~\ref{Lemma5QmRm} (x),(xi), $R_{M+1}\le Q_M<n+2.$ In either case, $n+1<R_{M+1}, Q_M<n+2$   and hence $n+1=\floor{Q_{M}}=\floor{R_{M+1}}.$

   \item Assume now $n-1\in [R_M,Q_M)$. \ If $n-1=R_{M}$ \ then \ $r=2n-1$  \, which, by Lemma~\ref{Lemma5QmRm}, is not possible as $M\ge 1.$ Hence $R_M<n-1.$  By Lemma~\ref{Lemma5QmRm} (ix), $n<Q_M.$ \ Therefore \ $n-2=\floor{Q_{M-1}}=\floor{R_M}.$ The situation looks like the one on Figure~\ref{FigureSecondSituation}

            \begin{figure}[H]
             \begin{center}
               \begin{tikzpicture}
                  \draw[thick] (0,0) -- (9,0);
                  \draw[thick] (4.5,.35) arc (160:200:1cm) node [below=5pt] {$R_M$};
                  \draw[thick] (2.5,.35) arc (20:-20:1cm) node [below=5pt] {$Q_{M-1}$};
                  \draw[thick] (7.5,.35) arc (20:-20:1cm) node [below=5pt] {$Q_M$};
                  \fill (1,0) circle (2.5pt) node [above=5pt] {$n-2$};
                  \fill (5,0) circle (2.5pt) node [above=5pt] {$n-1$};
                  \fill (6,0) circle (2.5pt) node [above=5pt] {$n$};
               \end{tikzpicture}
               \caption{The resulting case when $n-1\in [R_M,Q_M)$.}{\label{FigureSecondSituation}}
             \end{center}
            \end{figure}
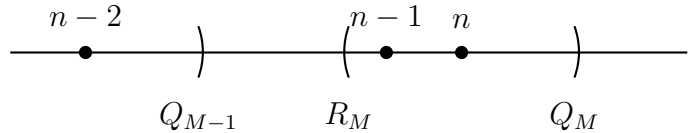

       Since $Q_{M-1}<n-1$, by Lemma~\ref{Lemma5QmRm} (viii),(x) \ $r>2n-3,\ Q_M<R_{M+1}.$ Since $R_M<n-1$, by Lemma~\ref{Lemma5QmRm} (ix), $R_{M+1}<n+1.$ Hence, $n<Q_M<R_{M+1}<n+1$   and hence $n=\floor{Q_{M}}=\floor{R_{M+1}}.$ \ok
 \end{itemize}


\section{Main Result}\label{Section:MainResult}

The following theorem provides a formula to compute the sum of the quadratic residues modulus $p=4n-1.$ Such a formula did not exist at the moment.
\begin{thm}\label{Theorem9Main4n-1} 
  Let $n\in\N, n>1$ such that $p=4n-1$ is prime. Using Notation~\ref{NotationQmRmM} we have

  \begin{equation}\label{FormulaSumofQRp=4n-1}
    \begin{split}
      \frac12\sum_{k=1}^{p-1}\r{k^2}
      \inviss{- p\sum_{m=0}^{M-1} (2m+1)\drqsp-2p\sum_{m=0}^{M-1}m\dqrsp}
             {+ \frac{4n^3+3n^2-n}{6}-Mp\Big(2n-1-2\fqMm\Big)}
    \end{split}
  \end{equation}
\end{thm}

\noindent {\bf Proof.}

  We observe that the sum of quadratic residues from 1 to $p-1$ is twice the sum of the quadratic residues from 1 to $2n-1$ as they repeat. Hence it is enough to look at the sum from 1 to $2n-1.$ \ Let $a_k=\r{k^2}+\r{(2n-k)^2}, \ b_k=2k^2-k+n$ \, and \, $c_k=a_k-b_k.$ In the sum of the quadratic residues from 1 to $2n-1$, if we add the first one and the last one, the second and the second to the last and so on, we obtain that
   \begin{equation*}
      \sum_{k=1}^{2n-1}\r{k^2}=a_n+\sum_{k=1}^{n-1} a_k=a_n+\sum_{k=1}^{n-1} b_k+\sum_{k=1}^{n-1} c_k.
   \end{equation*}

   From Corollary~\ref{Corollary3SquareRoots} we know that
   \[
      \defi{c_k}{-2mp & k\in(R_m,Q_m]\\ & \\ -(2m+1)p & k\in(Q_m,R_{m+1}).}
   \]

  From Lemma~\ref{Lemma5QmRm}, $n-1\in(R_{M-1},Q_M)-(Q_{M-1},R_M)=(R_{M-1},Q_{M-1}]\cup [R_M,Q_M).$  \ From the proof of Lemma~\ref{Lemma8EitherOr}, we know that $n-1\neq Q_{M-1}$ and $n-1\neq R_M,$ \, hence $n-1\in (R_{M-1},Q_{M-1})\cup (R_M,Q_M).$ \ Let \, $I=\{1,2,...,n-1\}.$
 \begin{itemize}
   \item  Case \, $n-1<Q_{M-1}.$

   By Lemma~\ref{Lemma8EitherOr}, $n-1=\floor{Q_{M-1}}=\floor{R_M}$ \, and also \, $R_{M-1}<n-2<n-1<Q_{M-1}.$
   Since all the $c_k$'s, for $k$ that are in certain intervals, have the same value, we have
   \begin{equation*}
    \begin{split}
      \sum_{k=1}^{n-1}c_k &=
               \sum\nolimits_{k\in I\cap\left(\bigcup_{m=0}^{M-1}(R_m,Q_m]\cap (Q_m,R_{m+1})\right)} c_k\\
      & =\sum_{m=0}^{M-1}\left(\sum_{k\in I\cap (R_m,Q_m]}      c_k
                                +\sum_{k\in I\cap (Q_m,R_{m+1})}  c_k\right)\\
      & =\sum_{m=0}^{M-1}\sum_{k\in I\cap (R_m,Q_m]} c_k
                                + \sum_{m=0}^{M-1}\sum_{k\in I\cap(Q_m,R_{m+1})} c_k\\
      & =\sum_{m=0}^{M-1}(-2mp)\hspace{-.7cm}\sum_{k\in I\cap (R_m,Q_m]}\hspace{-.5cm}1
                                + \sum_{m=0}^{M-1}(-(2m+1)p)\hspace{-.7cm}\sum_{k\in I\cap(Q_m,R_{m+1})} \hspace{-.7cm} 1\\
      & =-p\sum_{m=0}^{M-1} 2m\dqrsp
                                -p\sum_{m=0}^{M-1}(2m+1)\drqsp \\
      \inviss{ -p\sum_{m=0}^{M-1} 2m\dqrsp-p\sum_{m=0}^{M-1}(2m+1)\drqsp}{-2Mp\Big(n-1-\fqMm\Big).}\\
    \end{split}
   \end{equation*}

   \item Case \, $R_M<n-1.$  From the proof of Lemma~\ref{Lemma8EitherOr}, $\floor{Q_{M-1}}=\floor{R_M}=n-2$ \, $\floor{Q_{M}}=\floor{R_{M+1}}=n$ \, and also \, $R_{M-1}<n-2<Q_{M-1}<R_M<n-1<n<Q_M.$ \ Since $n-1\in (R_M,Q_M)$ then $a_{n-1}=2(n-1)^2-(n-1)+n-2Mp$, hence $c_{n-1}=-2Mp.$
   As in the previous case, to compute the sum of the $c_k$'s, we have the same terms plus this additional one. Note that $n-1-\fqMm=1.$
   \begin{equation*}
    \begin{split}
      \sum_{k=1}^{n-1}c_k &=
               \sum_{k=1}^{n-2}c_k + c_{n-1}\\
      \inviss{\sum_{m=0}^{M-1}(-(2m+1)p)\drq -2Mp}{+\sum_{m=0}^{M-1}(-2mp) \dqr}\\
      \inviss{-p\sum_{m=0}^{M-1}(2m+1)\drq -2Mp}{-p\sum_{m=0}^{M-1} 2m\dqr}\\
      \inviss{-p\sum_{m=0}^{M-1} 2m\dqr-2Mp\Big(n-1-\fqMm\Big)}{-p\sum_{m=0}^{M-1}(2m+1)\drq.}\\
    \end{split}
   \end{equation*}
 \end{itemize}

 Notice that $a_n=\frac12(\r{n^2}+\r{(2n-n)^2}=\frac12(2n^2-2Mp)=n^2-Mp.$ \ Since

   \begin{equation*}
    \sum_{k=1}^{n-1} b_k + n^2 =\sum_{k=1}^{n-1}(2k^2-k+n) +n^2=\frac{4n^3+3n^2-n}{6},
   \end{equation*}


we can now put all the pieces together to obtain,
\begin{equation*}
    \begin{split}
      \frac12\sum_{k=1}^{p-1}\r{k^2} & = \sum_{k=1}^{2n-1}\r{k^2}=a_n+\sum_{k=1}^{n-1} b_k+\sum_{k=1}^{n-1} c_k \\
         & = \frac{4n^3+3n^2-n}{6} -Mp+\sum_{k=1}^{n-1} c_k \\
     \inviss{-p\sum_{m=0}^{M-1} 2m\dqrsp-p\sum_{m=0}^{M-1}(2m+1)\drqsp}
            {+\frac{4n^3+3n^2-n}{6} -Mp -2Mp\Big(n-1-\fqMm\Big)}\\
     \inviss{- p\sum_{m=0}^{M-1} (2m+1)\drqsp-2p\sum_{m=0}^{M-1}m\dqrsp}
            {+ \frac{4n^3+3n^2-n}{6}-Mp\Big(2n-1-2\fqMm\Big)}.\\
    \end{split}
  \end{equation*} \ok

\begin{cor}\label{Corollary9.5TheoremMain4n-1} 
Under the hypothesis of Theorem~\ref{Theorem9Main4n-1},
  \begin{equation*}
     \frac12\sum_{k=1}^{p-1}\r{k^2} = p\left(\sum_{m=0}^{M} \frm +\sum_{m=0}^{M-1} \fqm\right)-Mp(2n-1)+\frac{p\cdot(n^2+n)}{6}
  \end{equation*}
\end{cor}
\noindent {\bf Proof.}
 From Theorem~\ref{Theorem9Main4n-1},
\begin{equation*}
    \begin{split}
      \frac12\sum_{k=1}^{p-1}\r{k^2}
      \inviss{- p\sum_{m=0}^{M-1} (2m+1)\drqsp-2p\sum_{m=0}^{M-1}m\dqrsp}
             {+ \frac{4n^3+3n^2-n}{6}-Mp\big(2n-1-2\fqMm\big)}\\
      \inviss{- 2p\sum_{m=0}^{M-1} m\drr-Mp\big(2n-1-2\fqMm\big)}
             {+\frac{4n^3+3n^2-n}{6}-p\sum_{m=0}^{M-2}\drq}\\
      \inviss{p\left(\sum_{m=0}^{M} \frm+\sum_{m=0}^{M-1} \fqm\right)-Mp\Big(2n-1-2\fqMm\Big)}
             {+\frac{4n^3+3n^2-n}{6}-2Mp\frM}.
    \end{split}
  \end{equation*}
  By Lemma~\ref{Lemma8EitherOr}, $\fqMm=\frM$, the result follows.
\ok

\section{A Computable Class Number Formula}\label{Section:Computable}
We now establish our formula that allows us to compute the class number for the imaginary quadratic field $\quadf$ when $p=4n-1$ is prime.

\begin{thm}\label{ComputableClassNumberFormula}
Under the hypothesis of Theorem~\ref{Theorem9Main4n-1},
  \[h(-p)=(2M+1)(2n-1) -2\left(\sum_{m=0}^{M} \frm +\sum_{m=0}^{M-1} \fqm\right)-\frac{(n^2+n)}{3}.\]
\end{thm}
\noindent {\bf Proof.}
According to the formula
\[\sum_{k=1}^{p-1} r_p(k^2)=\binom{p}{2}-p\cdot h,\]
and Corollary~\ref{Corollary9.5TheoremMain4n-1}, we obtain that
\[\binom{p}{2}-p\cdot h = 2p\left(\sum_{m=0}^{M} \frm +\sum_{m=0}^{M-1} \fqm\right)-2Mp(2n-1)+\frac{p\cdot(n^2+n)}{3},\]
the formula follows.
\ok

\begin{cor}\label{CorollaryClassNumberCases}
Under the hypothesis of Theorem~\ref{Theorem9Main4n-1}, if \[\gamma_k=\sum_{m=0}^{k} \frm +\sum_{m=0}^{k-1} \fqm \] then

  \[\defi{h(-p)}{
  \ds{\frac13} (32k^2+14k-3)-2\gamma_k, & p=16k-1\\[4mm]
  \ds{\frac13} (32k^2+18k+1)-2\gamma_k, & p=16k+3\\[4mm]
  \ds{\frac13} (32k^2+22k+3)-2\gamma_k, & p=16k+7\\[4mm]
  \ds{\frac13} (32k^2+26k+3)-2\gamma_k, & p=16k+11.} \]
\end{cor}
\noindent {\bf Proof.}
The proof is immediate by applying Theorem~\ref{ComputableClassNumberFormula} and Lemma~\ref{Lemma4.5Chart}.

\ok

A concise way to compute the class number is
\begin{cor}\label{EconomicCorollaryClassNumberCases}
Under the hypothesis of Theorem~\ref{Theorem9Main4n-1}, if \[\beta_0=2\floor{ \frac12+\frac12\newsqrt{3p-4}} \] then
  \[\defi{h(-p)}{
  \ds{\frac13\sum_{m=1}^{k}\Big(64m+6-6\frm-6\fqm\Big)-\beta_0+1}, & p=16k-1\\[4mm]
  \ds{\frac13\sum_{m=1}^{k}\Big(64m+10-6\frm-6\fqm\Big)-\beta_0+\frac{13}{3}}, & p=16k+3\\[4mm]
  \ds{\frac13\sum_{m=1}^{k}\Big(64m+14-6\frm-6\fqm\Big)-\beta_0+5}, & p=16k+7\\[4mm]
  \ds{\frac13\sum_{m=1}^{k}\Big(64m+18-6\frm-6\fqm\Big)-\beta_0+7}, & p=16k+11.} \]
\end{cor}
\noindent {\bf Proof.}
Observe that from Lemma~\ref{Lemma5QmRm} and Lemma~\ref{Lemma8EitherOr}, $r<2n-3$ implies $\fqM=n+1$ and  $r>2n-3$ implies $\fqM=n.$ Now, Table~\ref{TableResiduesn} in Lemma~\ref{Lemma4.5Chart} provides the different values for $r$ in each of the four cases and we notice that when $n=4k$ or $n=4k+1$, we have $r<2n-3$ and hence $\fqM=n+1.$ \ In the other two cases $r>2n-3$ and $\fqM=n.$ Consequently, if we look closely at Corollary~\ref{CorollaryClassNumberCases} and we write $\gamma_k$ and  $\fqM$ as well as the term containing $32k^2$  as a sum from 1 to $k,$ in we obtain  desired result. \ok


\section{Examples}\label{Section:Examples}
\begin{eg}\label{Example1}
As a simple example, consider $p=103=4\cdot 26-1$ a prime number, here $n=26=4\cdot 6+2.$
We compute\\[2mm]
    \begin{tabular}{|c|c|c|c|c|c|}
      \hline
      \sx $\fq{0}=9 $ & $\fq{1}=13 $ & $\fq{2}=17 $ & $\fq{3}=20 $ & $\fq{4}=22 $ & \sx $\fq{5}=24$ \\ \hline
      \sx $\fr{1}=10 $ & $\fr{2}=14 $ & $\fr{3}=17 $ & $\fr{4}=20 $ & $\fr{5}=22 $ & $\fr{6}=24$ \\
      \hline
    \end{tabular}\\[2mm]

 Hence by using Corollary~\ref{CorollaryClassNumberCases}, case $n=4k+2$, we obtain $k=6, \gamma_6=107+105=212$, and hence
 \[h(-103)=\ds{\frac13} (32k^2+22k+3)-2\gamma_k=429-424=5.\]
\end{eg}

\begin{eg}\label{Example1}
As a more complicated example, consider $k=50,000, \  n=4*k, p=4n-1=799,999$ a prime number.
We compute

\begin{eqnarray*}
  \beta_{0} &=& 1550, \\
  \sum_{m=1}^{50,000}\Big(64m+6-6\frm-6\fqm\Big) &=& 6216,
\end{eqnarray*}

\noindent and hence by using Corollary~\ref{EconomicCorollaryClassNumberCases}, case $n=4k$, we obtain
 \[h(-799,999)=\ds{\frac13}(6216)-1550+1=523.\]
 (We used Maple 2019 to compute the sum)
\end{eg}

\newpage
 \bibliographystyle{siam}
 \bibliography{quady}
\end{document}